\documentclass[8pt,a4paper]{article}
\usepackage{amsmath,amsfonts,amsthm,amsopn,color,amssymb,}
\newcommand{\eps}{\varepsilon}
\newcommand{\Rn}{{\mathbb{R}^n}}

\newcommand{\Rnu}{{\mathbb{R}^{n-1}}}

\renewcommand{\a }{\alpha }
\renewcommand{\le }{\leqslant }
\renewcommand{\ge }{\geqslant }
\renewcommand{\b }{\beta }
\renewcommand{\d }{\delta }

\renewcommand{\O}{\Omega}
\newcommand{\G}{\Gamma}

\renewcommand{\i}{{H}^{1}(\Omega)}

 \newcommand{\io}{\int _\O}
\newcommand{\iu}{\int _{\G_1}}
\newcommand{\id}{\int_{\G_2}}
\newcommand{\n}{\nabla}

\newtheorem{theorem}{Theorem}[section]
\newtheorem{lemma}[theorem]{Lemma}

\newtheorem{proposition}[theorem]{Proposition}

\newtheorem{corollary}[theorem]{Corollary}

\renewenvironment{proof}{\noindent{\textbf{Proof.\quad}}}{$\hfill\square$\vspace{0.2 cm}\\}
\newenvironment{proofb}{\noindent{\textbf{Proof of Theorem \ref{teo1}.\quad}}}{$\hfill\square$\vspace{0.2 cm}\\}
\newenvironment{proofa}{\noindent{\textbf{Proof of Theorem \ref{teo2}.\quad}}}{$\hfill\square$\vspace{0.2 cm}\\}

\begin{document}
{\title{ Detecting nonlinear corrosion \\ by electrostatic
measurements\thanks{Work supported in part by MIUR, grant n.
2002013279.} }}
\author{G. Alessandrini\thanks{Dipartimento di  Matematica e Informatica, Universita' degli Studi di Trieste, Italy,
\tt{alessang@univ.trieste.it}} \ and \ E.
Sincich\thanks{S.I.S.S.A., Trieste, Italy, \tt { sincich@sissa.it}
}}

\date{}

\maketitle

\begin{abstract}We deal with an inverse problem arising in corrosion detection. We
prove a stability estimate for a nonlinear term on the
inaccessible portion of the boundary by electrostatic boundary
measurements on the accessible one.
\end{abstract}



\section{Introduction}
In this paper we deal with the inverse problem of determining the
nonlinear term $f(u)$ in the boundary value problem
\begin{equation}\label{P}
\left\{
\begin{array}
{lcl}
\Delta u=0\ ,& \mbox{in $\Omega$ ,}
\\
\dfrac{\partial u}{\partial \nu}=g\ ,   & \mbox{on $\G_2$ ,}
\\
\dfrac{\partial u}{\partial \nu}=f(u)\ , & \mbox{on $\G_1$ ,}
\\
u=0\ , & \mbox{on $\G_{D}$ ,}
\end{array}
\right.
\end{equation}
where $\G_1$ and $\G_2$
are two open, disjoint portions of  $\partial\O$
and $\G_{D}=\partial \O \setminus (\G_1\cup\G_2)$ and $f$
is a Lipschitz
 function such that $f(0)=0$.

\noindent Such a problem with a specific choice of the profile of $f$, has been introduced and discussed in recent years by M.Vogelius and others for the modeling of the electrochemical phenomenon of surface corrosion in metals \cite{Vog1},\cite{Vog2},\cite{Vog3}.
 In the above boundary value problem, $\O$ represents the metal specimen, $\G_1$ represents the corroded part of the boundary, which is not accessible to direct inspection, $\G_2$ is the portion of the boundary which is accessible to measurements and $\G_D$ is the remaining part of $\partial \O$ which is assumed to be grounded.

\noindent The inverse problem thus consists in the determination of $f$ when one pair of Cauchy data $u|_{\G_2},{\frac{\partial u}{\partial \nu}}|_{\G_2}$ is known for one non-trivial solution $u$ to \eqref{P}.

\noindent It must be noted that the models of nonlinearity which have been discussed in \cite{Vog2},\cite{Vog3} are of the form
$$f(u)=\lambda(\exp(\a u)-\exp(-(1-\a)u))$$
and it turns out that are such that in the direct problem the
existence and the uniqueness of the solution are not granted. For
this reason, we found it necessary to require some additional a
priori assumptions. First an energy bound of the measured
electrostatic potential $u$
\begin{eqnarray}\label{E}
\int_{\O}|\nabla u(x)|^2\le E^2\ .
\end{eqnarray}
Next, an a priori bound of the Lipschitz continuity of $f$, namely
\begin{eqnarray}\label{Lip}
|f(u)-f(v)|\le L|u-v|\ , \ \mbox{for every}  \ u, v \ \in
\mathbb{R}\ .
\end{eqnarray}
Moreover, we found it necessary to assume the knowledge of some
additional information on the measured current density $g$ on the
accessible part of the boundary $\G_2$. More precisely, we assume
a bound on the H\"older continuity of $g$, with
\begin{eqnarray}\label{G}
 \|g\|_{C^{0,\a}(\G_2)}\le G\ .
\end{eqnarray}
Also, we shall require a lower bound on the same current density
$g$. Namely, we shall prescribe that, for a given inner portion
$\G_{2,2r_0}$ of $\G_2$ (see \eqref{gi} below, for a precise definition)
and a given number $m > 0$, we have
\begin{equation}\label{lower b}
\|g\|_{L^{\infty}(\G_{2,2r_0})}\ge m>0 \ .
\end{equation}

We note that, indeed, a lower bound of this sort appears to be
necessary. In fact, knowing the Cauchy data on $\G_2$ for a single
potential $u$, one can expect to identify $f$ only on the range of
values taken by $u$ on the inaccessible boundary $\G_1$. Thus, as
a preliminary step of the treatment of this inverse problem, it is
necessary to evaluate the amplitude of the range of $u$ on $\G_1$,
see Theorem \ref{teo1}, in which we prove that the oscillation of
$u$ on $\G_1$ is bounded from below by
$\exp(-\big(\frac{m}{c}\big)^{-\gamma})$, where $c>0,\gamma>1$ are
constants depending on the a priori data only. Next, as the main
result of this paper, we prove a stability estimate,  Theorem
\ref{teo2}, that is we show that if $u_1, u_2$ are two potentials
corresponding to nonlinearities $f_1,\ f_2,$ whose Cauchy data are
close
\begin{eqnarray*}
&&\|{\psi}_1-{\psi}_2\|_{L^2(\G_2)}\le \varepsilon\ ,\\
&&\|g_1-g_2\|_{L^2(\G_2)}\le \varepsilon\ ,
\end{eqnarray*}
where $\psi_i=u_1\big|_{\G_2}$ and $\displaystyle{\frac {\partial
u_i}{\partial \nu}}\bigg|_{\G_2}=g_i$ with $i=1,2$ and if $g=g_1$
satisfies the a priori bounds (\ref{G}), (\ref{lower b}), then the
ranges of $u_1, u_2$ on $\G_1$ agree on a sufficiently large
interval $V$ such that
$$\mbox{length of}\  V\  \sim \ \ \exp\bigg[- \bigg(\frac{m}{c}\bigg)^{-\gamma}\bigg] \ .$$
On such an interval the nonlinearities $f_1, f_2$ agree up to an
error $\omega(\varepsilon)$ such that
$$\omega(\varepsilon) \sim \left|\log\bigg(\frac{1}{\varepsilon}\bigg)\right|^{-\theta}$$
where $0<\theta<1$.

The study of determination of nonlinear terms in elliptic inverse
boundary value problems has produced various different results.
Just to mention a few, Cannon \cite{Ca}, Beretta and Vogelius
\cite{BV}, Isakov and Sylvester \cite{IS}, Sun \cite{Sun}. Among
the specific peculiarities of the present study, in comparison
with the mentioned previous results, we emphasize that the
nonlinearity is part of a boundary condition and that existence
and uniqueness for the corresponding direct problem are not
assured, see, in this respect, Kavian and Vogelius \cite{Vog2} and
Medville and Vogelius \cite{MV}. For the determination of a
coefficient in a \emph{linear} boundary condition, we refer for
instance to Fasino and Inglese \cite{FI}, Alessandrini, Del Piero
and Rondi \cite{adr} and Chaabane, Fellah, Jaoua and Leblond
\cite{CJ}.

\subsection{Plan of the paper}
In Section 2 we formulate our main hypotheses and state the main
results, Theorem \ref{teo1} and Theorem \ref{teo2}, whose proofs
are deferred to Section 3 and 4 respectively. In Section 3 we
prove some preliminary lemmas, in which we observe that, taking
into account the a priori assumptions, the solution is H\"older
continuous with its first order derivatives in a suitable
neighborhood of the inaccessible part $\G_1$, (see Lemma
\ref{localbound}, Theorem \ref{teoreg}). In order to prove Theorem
\ref{teo1} we need also a stability estimate near the boundary for
a Cauchy problem (see Proposition \ref{tritten}). The proof of
this Proposition is based on a method developed by Payne and
Trytten \cite{payne}, \cite{tr}. We conclude the Section with the
proof of Theorem \ref{teo1}. In Section 4 we adapt the above
mentioned stability estimate for the Cauchy problem (see Theorem
\ref{stabilita'}), when both the Dirichlet and Neumann data are
affected by errors in the $L^2$ norm, and at the same time, we
obtain a stability result up to the boundary. In Proposition
\ref{inv} we show that the solution is locally invertible along a
suitable curve on $\G_1$ and we evaluate its length. Finally we
prove Theorem \ref{teo2}.


\section{Main assumptions and results}
\subsection{Preliminaries}
{\it A priori information on the domain}
\\
\\
We shall assume throughout that $\O$ is a bounded, simply
connected domain in $\Rn$, $n\ge2$ such that $\mbox{diam} (\O)\le
D$ and with Lipschitz boundary $\partial \O$ with constants $r_0,
M$. More precisely, for every $x_0 \in \partial \O$, exists a
rigid transformation of coordinates under which,
\begin{eqnarray}\label{gamma}
\O \cap B_{r_0}(x_0)=\{(x',x_n): x_n>\gamma(x')\}
\end{eqnarray}
where $x \in \Rn,\ x=(x',x_n)$, with $x'  \in \Rnu , \ x_n \in \mathbb{R}$ and
$$\gamma :B_{r_0}(x_0)\subset \Rnu \rightarrow \mathbb{R}\ ,$$
satisfying $\gamma(0)=0$ and $$\|\gamma\|_{C^{0,1}( B_{r_0}(x_0))}\le Mr_0\ ,$$ where we denote by
$$\|\gamma\|_{C^{0,1}( B_{r_0}(x_0))}=\|\gamma\|_{L^{\infty}( B_{r_0}(x_0))} +\ r_0\!\!\!\!\!\!\!\!\!\sup_{\substack {x,y  \in B_{r_0}(z_0)\\x\ne y }}
\frac{|\gamma (x)-\gamma (y)|}{|x-y|}\ . $$
Moreover, we assume that the portions of the boundary $\Gamma_i$ are contained respectively into surfaces $S_i$, $i=1,2$ which are $C^{1,\a}$ smooth with constants $r_0, M$.

More precisely,
 for a given $\a,\
 0<\a<1$ and for any $z_0 \in S_i,\  i=1,2 $, we have that up to a rigid change of coordinates,
\begin{eqnarray}\label{si}
 S_i \cap B_{r_0}(z_0)=\{(x',x_n): x_n=\varphi_i(x')\}\ ,
\end{eqnarray}
where
$$\varphi_i:B_{r_0}(z_0)\subset \Rnu \rightarrow \mathbb{R} $$
 are  $C^{1,\a}$ functions satisfying $\varphi_i(0)=|\nabla \varphi_i(0)|=0$ and $$\|\varphi_i \|_{C^{1,\a}(B_{r_0}(z_0))}\le Mr_0\ ,$$
where we denote
\begin{eqnarray}\nonumber
\|\varphi\|_{C^{1, \a}( B_{r_0}(z_0))}&=&\|\varphi\|_{L^{\infty}( B_{r_0}(z_0))}+r_0\|\nabla \varphi\|_{L^{\infty}( B_{r_0}(z_0))}+\\
& & +\ {r_0}^{1+\a}\!\!\!\!\!\!\!\!\!\sup_{\substack {x,y  \in B_{r_0}(z_0)\\x\ne y }}\frac{|\nabla \varphi (x)-\nabla \varphi (y)|}{|x-y|^{\a}}\ .
\end{eqnarray}

In particular it follows that if
$$x_0 \in \Gamma_i\ \ \mbox{and} \ \ \mbox{dist}(x_0,\G_D)>r_0\ ,$$
then
$$\O \cap B_{r_0}(z_0)=\{(x',x_n) \in  B_{r_0}(z_0): x_n>\varphi_i(x')\}\ ,$$
where $\varphi_i$ is the Lipschitz function whose graph locally represents $\partial \O$. Moreover, since $\O \cap B_{r_0}(z_0) \cap \G_D =\emptyset$, $\varphi_i$ must also be the $C^{1,\a}$ function whose graph locally represents $S_i$.
We also suppose that the boundary of $\G_i$, within $S_i$, is of $C^{1,\a}$ class with constants $r_0, M$, namely, for any $z_0\in \partial \G_i$, there exists a rigid transformation of coordinates under which
$$\partial \G_i \cap B_{r_0}(z_0)=\{(x',x_n)\in B_{r_0}(z_0):\ \ x_n=\varphi_i(x'),\  x_{n-1}=\psi_i(x'')\}$$
where $x'=(x'',x_{n-1})$, with $x'' \in \mathbb{R}^{n-2}$, $x_{n-1} \in \mathbb{R}$ and
$$\psi_i : B_{r_0}(z_0)\subset  \mathbb{R}^{n-2} \longrightarrow \mathbb{R} $$
satisfying $\psi_i(0)=|\nabla \psi_i(0)|=0$ and $\|\psi_i\|_{C^{1,\a} (B_{r_0}(z_0))}\le M$.

We introduce some notation that we shall use in the sequel
\begin{eqnarray}\label{ui}
U^i_{\rho}=\{x \in \bar{\O}:\mbox{dist}(x,\partial\O\setminus \G_i)>\rho \}\ \ ,
\end{eqnarray}
\begin{eqnarray}\label{gi}
\Gamma_{i,\rho}=U^1_{\rho}\cap \G_i\ ,
\end{eqnarray}
with i=1,2  and
\begin{eqnarray}\label{or}
\O_{\rho}=\{x\in \O: \mbox{dist}(x, \partial \O)>\rho \}\ .
\end{eqnarray}

 {\it A priori information on the boundary data}
\\
\\
The current flux $g$
is a prescribed function such that
\begin{eqnarray*}
 \|g\|_{C^{0,\a}(\G_2)}\le G\ .
\end{eqnarray*}

 {\it A priori information on the nonlinear term}
\\
\\
We assume that the function $f$ belongs to $C^{0,1}(\mathbb{R},\mathbb{R})$ and, in particular,
\begin{eqnarray}\label{f}
  f(0)=0\ \  \mbox{and}\ \  | f(u)-f(v)|\le L|u-v|\ \ \
 \mbox{for every } u, v\in \mathbb{R}\ .
\end{eqnarray}

We recall that a weak solution of problem \eqref{P} is a  function $ u \in \i$, such that $u|_{\G_D}=0$ and which satisfies
\begin{eqnarray}\label{soluzione}
\io \n u\cdot \n \rho  =\id g  \rho + \iu f(u) \rho \quad \mbox{for all } \rho \in \i, \, \,\rho _{|_{_{_{\G_D}}}}=0
\end{eqnarray}

where  $u_{|_{_{_{\G_D}}}}  $ denote the trace on $\G_D$ .
\\
From now on we shall refer to the \emph{a priori data} as to the set of quantities
$r_0, M, \a, L, G, E, D, m$.\\
In the sequel we shall denote with $\eta(t)$ and $\omega(t)$, two positive increasing functions defined on $(0, +\infty)$, that satisfy
\begin{eqnarray}\label{eta}
&&\eta(t)\ge\exp{\bigg[-\bigg(\frac{t}{c}\bigg)^{-\gamma}\bigg]}, \ \ \ \mbox{for every}\ \ 0<t\le G\ \ ,
\end{eqnarray}
\begin{eqnarray}\label{omega}
&&\omega(t)\le C \left|\log (t) \right|^{-\vartheta},\ \ \ \mbox{for every}\ \ 0<t<1\ \ ,
\end{eqnarray}
where $c>0$, $C>0$, $\gamma>1$, $0<\theta<1$ are constants depending on the \emph{a priori data} only.
\subsection{Main theorems}
The statements of the main results are the following.
\begin{theorem}[Lower bound for the oscillation]\label{teo1}
Let $\O, g $ satisfying the a priori assumptions. Let $u$ be a
weak solution of \eqref{P} satisfying the a priori bound \eqref{E}
then $$\mathop{\rm osc}\limits_{\G_1} u\ge
\eta(\|g\|_{L^{\infty}(\G_{2,2r_0})})$$ where $\eta$ satisfies
\eqref{eta}.
\end{theorem}
\begin{theorem}[Stability for the nonlinear term $f$]\label{teo2}
Let $u_i\ \in \i$, $i=1,2$ be weak solutions of the problem $\eqref{P}$, with $f=f_i$ and $g=g_i$
respectively and such that \eqref{E} holds for each $u_i$.
Let us also assume that, for some positive number $m$, the following holds
\begin{eqnarray}\label{m}
\|g_1\|_{L^{\infty}(\G_{2,2r_0})}\ge m>0\ .
\end{eqnarray}
Moreover, let $\psi_i=u_i\big|_{\G_2}, \ i=1,2$.
There exists $\eps_{0}>0$ only depending on the \emph{a priori data} and on $m$ such that, if, for some $\eps,\ 0<\eps<\eps_0$, we have
\begin{eqnarray*}
&& \|\psi_1 -\psi_2  \|_{L^{2}(\G_2)}\le \eps\ , \\
&& \|g_1 -g_2  \|_{L^{2}(\G_2)}\le \eps\ ,
\end{eqnarray*}
 then
$$\|f_1(u)-f_2(u) \|_{L^{\infty}(V)}\le \omega(\eps)\ , $$
where
 $$V=(\a,\beta)\subseteq [-CE,CE]\ ,$$
is  such that
 $$\beta -\a >\frac{\eta(m)}{2}$$
and $\eta, \omega$ satisfy \eqref{eta}, \eqref{omega} respectively.
\end{theorem}

\section{The lower bound for the oscillation}

\begin{proposition}[Stability for a Cauchy problem]\label{tritten}
Let $\O$ satisfies the a priori assumptions and let $v \in
H^1(\O)$ be a solution of the following Cauchy problem
\begin{equation}
\left\{
\begin{array}{lcl}
\Delta v=0\ ,& \mbox{in $\Omega$},
\\
v=\varphi\ ,  & \mbox{on $\G_1$},
\\
\displaystyle\frac{\partial v}{\partial \nu}=h, & \mbox{on $\G_1$},
\end{array}
\right.
\end{equation}
where $\varphi, h\in L^2(\G_1)$ and the boundary conditions are considered in the weak sense.\\
Then, for every $P_1 \in \G_{1,2r_0}$, $v$ satisfies the following estimate
\begin{eqnarray*}
\|v\|_{L^2(B_{\rho}(P_0)\cap U^1_{2r_0})}&\leq&\ C \left(\|\varphi\|_{L^2(\G_{1,\rho})} +\
\|h\|_{L^2(\G_{1,\rho})}\ +\|v\|_{\i}
\right)^{1-\delta}\cdot\\
&&\ \ \ \ \ \ \ \ \ \cdot\left(\|\varphi\|_{L^2(\G_{1,\rho})} +\
\|h\|_{L^2(\G_{1,\rho})}\right)^{\delta}
\end{eqnarray*}
where $ \rho \in \left(\frac{M}{4\sqrt{1+M^2}}r_0,\frac{3M}{4\sqrt{1+M^2}}r_0 \right)$, $P_0=P_1+\frac{M}{4\sqrt{1+M^2}}r_0\cdot \nu$, $\nu$ is the outer unit normal to $\O$ at $P_1$ and
 $C>0$, $0<\delta<1$ are constants depending on
 $ \rho,r_0, n, M $ only.
\end{proposition}

\begin{proof}By a well-known estimate of stability for a Cauchy problem, we have that for any $P_1 \in \G_{1,r_0}$
\begin{eqnarray}\label{grad}
\|v\|_{L^2(B_{\rho}(P_0)\cap U^1_{r_0})}&\leq&C_1 \left(\|\varphi\|_{L^2(\G_{1,\rho})} +\
\|\nabla v\|_{L^2(\G_{1,\rho})}\ +\|v\|_{\i}
\right)^{1-\delta}\cdot \nonumber\\
&&\ \ \ \ \ \ \ \ \ \cdot\left(\|\varphi\|_{L^2(\G_{1,\rho})} +\
\|\nabla v\|_{L^2(\G_{1,\rho})}\right)^{\delta}
\end{eqnarray}
where $\rho \in \left(\frac{M}{4\sqrt{1+M^2}}r_0,\frac{3M}{4\sqrt{1+M^2}}r_0 \right)$ and $P_0=P_1+\frac{M}{4\sqrt{1+M^2}}r_0\cdot \nu$, where $\nu$ is the outer unit normal to $\O$ at $P_1$ and $C_1>0$, $0<\eta<1$ are constants depending on
 $ \rho,r_0, n, M $ only.
For a  proof we refer to \cite{tr}, see also \cite{payne}.
Let us define the function $\tilde{h} \in L^2(\partial \O)$ as follows
$$
\tilde{h}(x)=
\left\{
\begin{array}{ll}
h(x), &\quad\mbox{for a.e. $x \in \G_{1,r_0}$}\ , \\
\displaystyle-\frac{1}{|\G_2|}\int_{\G_{1,r_0}}h\ , & \quad\mbox{for a.e. $x \in \G_2$}\ ,\\0\ , & \quad\mbox{otherwise .}
\end{array}
\right.
$$
Let us consider the following Neumann problem
\begin{equation}\label{z}
 \left\{
\begin{array}
{lcl}
\Delta z=0\ ,& \mbox{in $\Omega$}\ ,
\\
\dfrac{\partial z}{\partial \nu}=\tilde{h}\ , & \mbox{on $\partial \O$ .}
\end{array}
\right.
\end{equation}
Note that $\int_{\partial \O}\tilde{h}=0$, hence a weak solution $z\in H^1(\O)$ exists and it is unique up to an additive constant.
We select the solution $z$ of \eqref{z} with zero average, it is well-known that for such a $z$ the following holds
$$\|z\|_{\i}\le C_2\|\tilde{h}\|_{L^2(\partial \O)}\le C_3\|h\|_{L^2(\G_{1,r_0})}$$
where $C_2$ and $C_3$ are positive constants depending on the \emph{a priori data} only.
Let us set $w=v-z$, thus $w$ solves the following Cauchy problem
\begin{equation}\label{w}
 \left\{
\begin{array}
{lcl}
\Delta w=0\ ,& \mbox{in $\Omega$},
\\
w=\varphi-z\ , & \mbox{on $\G_{1,r_0}$},
\\
\dfrac{\partial w}{\partial \nu}=0\ , & \mbox{on $\G_{1,r_0}$.}
\end{array}
\right.
\end{equation}
By a standard boundary regularity estimate (see for instance
\cite[p.667] {adn}), we have that $w \in
C^{1,\beta}(U^1_{\frac{3}{2}r_0})$ and the following holds
\begin{eqnarray}\label{beta}
\|w\|_{C^{1,\beta}(U^1_{\frac{3}{2}r_0})}\le C_4\|w\|_{\i}
\end{eqnarray}
where  $0<\beta<1$ and $C_4>0$ are constants depending on
$r_0,M,\a$ only. By an interpolation inequality, (see for instance
\cite[p.777]{ABRV}) we have that
\begin{eqnarray}\label{interpol}
\|\nabla w\|_{{L^2}(\G_{1,2r_0})}\le C_5\|w\|^{1-\gamma}_{C^{1,\beta}(U^1_{\frac{3}{2}r_0})}\|w\|^{\gamma}_{L^2(\G_{1,2r_0})}
\end{eqnarray}

where $C_5>0$ and $0<\gamma<1$ are constants depending on $M,\a,r_0$ only.

Moreover,
\begin{eqnarray}\label{distriang}
\|w\|_{\i}\le \|v\|_{\i}+\|z\|_{\i}\le C_6\left(\|v\|_{\i}+\|h\|_{L^2(\G_{1,r_0})}\right)
\end{eqnarray}
where $C_6=\max\{1,C_3 \}$.
From \eqref{beta},\eqref{interpol} and \eqref{distriang} it follows that
\begin{eqnarray}\label{nonso}
\|\nabla w\|_{L^{2}(\G_{1,2r_0})}\le C_7\left(\|v\|_{\i}+\|h\|_{L^2(\G_{1,r_0)})}\right)^{1-\gamma}\cdot\\\left(\|\varphi\|_{L^2(\G_{1,2r_0})}+\|z\|_{L^2(\G_{1,2r_0})} \right)^{\gamma}\ .
\end{eqnarray}

Applying \eqref{grad} to $w$ and using \eqref{nonso} we obtain
\begin{eqnarray*}
\|v\|_{L^2(B_{\rho}(P_0)\cap U^1_{2r_0})}&\leq&C \left(\|\varphi\|_{L^2(\G_{1,r_0})} +\
\|h\|_{L^2(\G_{1,r_0})}\ +\|v\|_{\i}
\right)^{1-\gamma\eta}\cdot\\
&&\cdot\left(\|\varphi\|_{L^2(\G_{1,r_0})} +\
\|h\|_{L^2(\G_{1,r_0})}\right)^{\gamma\eta}\ .
\end{eqnarray*}
And the thesis follows with $\delta=\gamma\eta$.
\end{proof}

\begin{lemma}[H\"older regularity at the boundary]\label{localbound}
Let $u$ be a solution of \eqref{P}, satisfying the  a priori bound \eqref{E}
then there exists a constant $C>0$, depending on the \emph{a priori data} only, such that
\begin{eqnarray}\label{x}
 \|u\|_{L^{\infty}(B_{\frac{r_0}{4}}(z_0)\cap \O)}\le CE ,\  \  \mbox{for every }  \ z_0\in \G_1
\end{eqnarray}
and
\begin{eqnarray}\label{y}
\|u\|_{C^{0,\a}(\G_1)}\le C E
\end{eqnarray}
where $0<\a<1$ is a constant depending on $r_0,M,n$ only.
\end{lemma}

\begin{proof}
From the weak formulation of the problem \eqref{soluzione} we have
\begin{eqnarray*}
\int_{B_{\frac{r_0}{2}}}\nabla u\cdot \nabla \varphi=\int_{\G_{\frac{r_0}{2}}}f(u)\varphi
\end{eqnarray*}
where $B_{\frac{r_0}{2}}=B_{\frac{r_0}{2}}(z_0) \cap \O$ and $\G_{\frac{r_0}{2}}=\partial
B_{\frac{r_0}{2}} \cap \G_1$ and $\varphi$ is any test function in $\i$ such that $\mbox{supp}{\varphi}\subset B_{\frac{r_0}{2}} \cup \G_{\frac{r_0}{2}}$.\\
By \eqref{f} we have that
\begin{eqnarray*}
\bigg|\int_{B_{\frac{r_0}{2}}}\nabla u\cdot \nabla \varphi\bigg|\le L\int_{\G_{\frac{r_0}{2}}}|u\varphi|
\end{eqnarray*}
and by a trace inequality (see \cite[p.114]{Adams}) it follows that
\begin{eqnarray}\label{traccia}
\bigg|\int_{B_{\frac{r_0}{2}}}\nabla u\cdot \nabla \varphi\bigg|\le CL \int_{B_{\frac{r_0}{2}}}|\nabla (u\varphi)|.
\end{eqnarray}
From now on one can proceed by the standard iteration technique due to Moser (see for instance \cite {gt}) and \eqref{x} follows. By the local bound \eqref{x} and by applying again the method by Moser leading to the Harnack inequality (see \cite{gt}) we obtain \eqref{y}.
 \end{proof}

\begin{theorem}[$C^{1,\a}$ regularity at the boundary]\label{teoreg}
Let $u$ be a solution of \eqref{P}, satisfying the  a priori bound \eqref{E}, then $u \in C^{1,\alpha}(G)$ and there exists a constant $C_{\rho}>0$, depending on the \emph{a priori data} and on $\rho$ only, such that the following estimate holds
\begin{eqnarray}\label{stima2}
\|u \|_{C^{1,\alpha}(U^1_{\rho})}\le C_{\rho}E
\end{eqnarray}
where  $\rho \in (0,r_0)$ .
\end{theorem}
\begin{proof}
 Since, by Lemma \ref{localbound}, we know that $u \in  C^{0,\alpha}(\G_1)$, by the Lipschitz
 regularity of $f$ we have that
$$\frac{\partial u}{\partial \nu}(x)=f(u(x))\ \ \in  C^{0,\alpha}(\G_1)\ . $$
By well-known regularity bounds for the Neumann problem (see for instance \cite[p.667]{adn}) it follows that $u \in C^{1,\alpha}(U^1_{\rho})$ and the following estimate holds
\begin{eqnarray}\label{stimaadn}
\|u\|_{C^{1,\alpha}(U^1_{\rho})}&\le& C\left(\|u\|_{C^{0,\alpha}(\G_{1,\frac{\rho}{2}})}+\left\|\frac{\partial u}{\partial \nu}\right\|_{C^{0,\alpha}(\G_{1,\frac{\rho}{2}})}+\|\nabla u\|_{L^2(\O)}\right)\le\nonumber\\
&\le&C\left(\left\|\frac{\partial u}{\partial \nu}\right\|_{C^{0,\alpha}(\G_{1,\frac{\rho}{2}})}+E \right)
\end{eqnarray}
where $C>0$ depends on the \emph{a priori data} and on $\rho$ only.
Moreover, we can estimate the $C^{0,\alpha}$ norm of $\displaystyle\frac{\partial u}{\partial \nu}$ in terms of $E$, in fact

\begin{eqnarray}
\left\|\frac{\partial u}{\partial \nu}\right\|_{C^{0,\alpha}(\G_{1,\frac{\rho}{2}})}=
\sup_{x \in \G_{1,\frac{\rho}{2}}}\left|\frac{\partial u(x)}{\partial \nu}\right|+\left({\frac{\rho}{2}}\right)^{\a}\sup_{x,y \in \G_{1,\frac{\rho}{2}}}\frac{\left|\frac{\partial u(x)}{\partial \nu}-\frac{\partial u(y)}{\partial \nu}\right|}{|x-y|^{\alpha}}=\nonumber \\
=\sup_{x \in \G_{1,\frac{\rho}{2}}}\left|f(u(x))\right|+\left({\frac{\rho}{2}}\right)^{\a}\sup_{x,y \in \G_{1,\frac{\rho}{2}}}\frac{\left|f(u(x))-f(u(y))\right|}{|x-y|^{\alpha}}\ .\nonumber
\end{eqnarray}
 By the Lipschitz bound \eqref{Lip} on $f$ and by Lemma \ref{localbound} we obtain
\begin{eqnarray}
\left\|\frac{\partial u}{\partial \nu}\right\|_{C^{0,\alpha}(\G_{1,\frac{\rho}{2}})}&\le& L\sup_{x \in \G_{1,\frac{\rho}{2}}}|u(x)|+L\left({\frac{\rho}{2}}\right)^{\a}\sup_{x,y \in \G_{1,\frac{\rho}{2}}}\frac{|u(x)-u(y)|}{{|x-y|}^{\alpha}}\le\nonumber\\
&\le& L\|u \|_{C^{0,\alpha}(\G_1)}\le CE\ .
\end{eqnarray}
So inserting this estimate in \eqref{stimaadn} we have the thesis.
\end{proof}
\begin{corollary}\label{corollario}
Let $u$ be as above, then, for every $\rho>0$, the function $\frac{\partial u}{\partial \nu}$ belongs to $C^{0,1}(\G_{1,\rho})$, with Lipschitz constant $\tilde{L}$ depending on the a \emph{a priori data} and on $\rho$ only.
\end{corollary}
\begin{proof}
Let $x$ and $y$ be two points in $\G_{1,\rho}$ then, by the assumption \eqref{Lip} and by Theorem \ref{teoreg}, it follows that
\begin{eqnarray*}
\bigg|\frac{\partial u(x)}{\partial \nu}-\frac{\partial u(y)}{\partial \nu}\bigg|&=& |f(u(x))-f(u(y))|\le L |u(x)-u(y)|\le\\
&\le& LC_{\rho}E|x-y|\ .
\end{eqnarray*}
The thesis follows with $\tilde{L}= LC_{\rho}E$.
 \end{proof}

\begin{proofb}
Let  $\eps=\mathop{\rm osc\ u}\limits_{\G_1}>0 $, since $u=0$ on $\G_D$, we have that
\begin{eqnarray}\label{bu}
 \|u\|_{L^2({\G_{1,r_0})}}\le C_1\eps
\end{eqnarray}
where $C_1$ is a positive constant depending on the \emph{a priori
data} only. By the a priori assumption \eqref{f} on $f$, we have
that  $|f(u)|\le L\eps $, moreover, since $$\left|\frac{\partial
u(x)}{\partial \nu}\right|=|f(u(x))|\ \  \mbox{on   $\G_1$},$$
then
\begin{eqnarray}\label{bu2}
\left\|\frac{\partial u}{\partial \nu}\right \|_{L^2(\G_{1,r_0})}\le |\G_{1,r_0}|^{\frac{1}{2}}L\eps\ .
\end{eqnarray}
By Proposition \ref{tritten}, it follows
\begin{equation}\label {tri}
 \|u\|_{L^2(B_{\rho}(P_0)\cap U^1_{2r_0})}\le C (\eps+E)^{1-\d}\cdot \eps ^{\d}
\end{equation}
where $C$ is a constant depending on the \emph{a priori data} only.
Since the boundary of $\O$ is of Lipschitz class, then it satisfies the cone property.
More precisely, if $Q$ is a point of $\partial \O$, then there exists a rigid transformation of coordinates under which we have $Q=0$. Moreover, considering the finite cone
$$\mathcal{C}=\bigg\{x :|x|<r_0 ,\,\,\frac{x\cdot \xi}{|x|}>\cos\theta \bigg \} $$
with axis in the direction $\xi$ and width $2\theta$, where $\theta=\arctan \frac {1}{M}$, we have that $\mathcal{C}\subset \O$. Let us consider now a point $Q \in \G_{2,r_0}$  and let $Q_0$ be a point lying on the axis $\xi$ of the cone with vertex in $Q=0$ such that $d_0=\mbox{dist}(Q_0,0)<\frac{r_0}{2}$. Following Lieberman \cite{Li}, we introduce a regularized distance $\tilde{d}$ from the boundary of $\O$. We have that there exists $\tilde{d}$ such that $\tilde{d}\in C^2(\O)\cap C^{0,1}({\bar\O})$, satisfying the following properties
\begin{itemize}
\item $\gamma_0\le \displaystyle \frac{\mbox{dist}(x,\partial \O)}{\tilde{d}(x)}\le \gamma_1$,
\item $|\nabla {\tilde{d}}(x)|\ge c_1$,\ \ \ for every $x$ such that ${\mbox {dist}}(x,\partial \O)\le br_0$,
\item $\|\tilde{d}\|_{C^{1,\a}}\le c_2r_0 $,
\end{itemize}
where $\gamma_0, \gamma_1, c_1, c_2, b$ are positive constants depending on $M, \a$ only, (see also \cite[Lemma 5.2]{ABRV}).

Let us define for every $\rho>0$
$$\tilde\O_{\rho}=\{x\ \in \O\ \ : \ \ {\tilde{d}}(x)>\rho \}\ . $$
It follows that, there exists $a$, $0<a\le1$, only depending on $M, \a$ such that for every $\rho$, $0<\rho\le ar_0$, $\tilde\O_{\rho}$ is connected with boundary of class $C^1$
and
\begin{eqnarray}\label{dc}
{\tilde c_1}\rho\le \mbox{dist}(x,\partial\O)\le \tilde c_2\rho \ \ \ \ \mbox{for every}\ \ x \ \in \ \partial\tilde\O_{\rho} \cap \O
\end{eqnarray}
where $\tilde c_1, \tilde c_2$ are positive constants depending on $M,\a$ only.
By \eqref{dc} it follows that
$$\O_{\tilde c_2\rho}\subset \tilde \O_{\rho}\subset \O_{\tilde c_1\rho}\ . $$
Using the notation introduced in the Proposition \ref{tritten}, we define the point $P=P_0-\frac{1}{4\sqrt{1+M^2}}r_0\cdot\nu$ and  ${\rho}_0=\min\{\frac{1}{32M\sqrt{1+M^2}}r_0,\frac{r_0}{4}\sin \theta\}$.
Moreover, let  $\gamma$ be a path in $\tilde\O_{\frac{\rho_0}{\tilde c_1}} $ joining  $P$ to  $Q_0$ and let us define $\{y_i\}$, $ i=0,\ldots,s$
as follows $y_0=Q_0$, $y_{i+1}=\gamma(t_{i})$, where $t_i=\max\{\mbox{$t$ s.t. }|\gamma(t)-y_i|=2\rho_0\}$
if $|P-y_i|>2\rho_0$ otherwise let $i=s$
and stop the process.

Now, we will use the three spheres inequality for harmonic functions (see for instance \cite{KM} or \cite[Appendix E]{ab}) that is
$$ \int_{B_{3\rho_0}(y_0)} u^2 \le \left(\int_{B_{\rho_0}(y_0)} u^2\right)^{\tau}\cdot \left(\int_{B_{4\rho_0}(y_0)} u^2\right)^{1-\tau}$$

where $0<\tau<1$ is an absolute constant. Now since $B_{\rho_0}(y_0)\subset B_{3\rho_0}(y_1)$ and since, by hypothesis  $\|u\|_{\i} \le E$, then we have

$$ \int_{B_{\rho_0}(y_0)} u^2 \le \left(\int_{B_{3\rho_0}(y_1)} u^2\right)^{\tau}\cdot E^{1-\tau}\ .$$
An iterated application of the three spheres inequality leads to
$$  \int_{B_{\rho_0}(y_0)} u^2 \le   \left(\int_{B_{\rho_0}(y_s)} u^2\right)^{{\tau}^{s}}\cdot E^{1-{\tau}^s}\ .$$
Finally, since we have $B_{\rho_0}(y_s)\subset B_{\frac{3M}{4\sqrt{1+M^2}}r_0}(P_0)\cap U^1_{2r_0}$, then by the Proposition \ref{tritten} it follows
$$ \int_{B_{\rho_0}(y_0)} u^2 \le C\big\{ (\eps+E)^{1-\d}\cdot(\eps )^{\d}\big\}^{{\tau}^s}\ .$$
We shall construct a chain of balls $B_{\rho_k}(Q_k)$ centered on the axis of the cone, pairwise tangent to each other and all contained in the cone
$${\mathcal{C}}^{\prime}=\bigg\{x :|x|<r_0 ,\,\,\frac{x\cdot \xi}{|x|}>\cos{\theta}^{\prime} \bigg \} $$
where $\theta^{\prime}=\arcsin\big(\frac{\rho_0}{d_0}\big) .$
Let $B_{\rho_0}(Q_0)$ be the first of them, the following are defined by induction in such a way
\begin{equation*}
\begin{array}{l}
Q_{k+1}=Q_{k}-(1+\mu)\rho_k\xi\ ,
\\
\rho_{k+1}=\mu\rho_{k}\ ,
\\
d_{k+1}=\mu d_{k}\ ,
\end{array}
\end{equation*}
with
\begin{equation*}
\begin{array}{l}
 \mu=\dfrac{1-\sin\theta^{\prime}}{1+\sin\theta^{\prime}}\ .
\end{array}
\end{equation*}
Hence, with this choice, we have $\rho_{k}=\mu^{k}\rho_0$ and  $B_{\rho_{k+1}}(Q_{k+1})\subset B_{3\rho_k}(Q_k)$.\\
Let us now consider the following estimate obtained by a repeated application of the three spheres inequality
\begin{eqnarray}\label{s}
\| u \|_{L^2(B_{\rho_k}(Q_k))}&\le& \| u \|_{L^2(B_{3\rho_{k-1}}(Q_{k-1}))
}\le \nonumber\\
&\le&  \| u \|^{\tau}_{L^2(B_{\rho_{k-1}}(Q_{k-1}))}  \| u \|^{1-\tau}_{L^2(B_{4\rho_{l-1}}(Q_{k-1}))}\nonumber\\
&\le& C \| u  \|^{{\tau}^k}_{L^2(B_{{\rho}_0}(Q_0))}\le\nonumber\\
&\le& C\Big\{ \big[ (\eps+E)^{1-\d}\cdot(\eps )^{\d}\big]^{{\tau}^s}\Big\}^{{\tau}^k}\ .
\end{eqnarray}
For every $r$, $0<r<d_0$, let $k(r)$ be the smallest positive integer such that  $d_k\le r$, then since $d_k={\mu}^k d_0$, it follows
\begin{eqnarray}\label{r}
\dfrac{|\log(\frac{r}{d_0})|}{\log{\mu}}\le k(r)\le \dfrac{|\log(\frac{r}{d_0})|}{\log{\mu}} +1
\end{eqnarray}
and by $(\ref{s})$, we have
\begin{eqnarray}\label{k}
   \| u \|_{L^2(B_{\rho_k(r)}(Q_k(r)))}\le
  C  \Big\{ \big[ (\eps+E)^{1-\d}\cdot(\eps)^{\d}\big]^{{\tau}^s}\Big\}^
{{\tau}^{k(r)}}\ .
\end{eqnarray}
Since, by hypothesis, $\G_2$ is contained in a $C^{1,\a}$ surface and by the regularity assumption \eqref{G} on $g$, it follows, by the same argument used in Theorem \ref{teoreg}, that  $u \in  C^{1,\a}(U^2_{2r_0}) $.

Let $\bar{x} \in \G_{2,2r_0}$, $x \in B_{\frac{\rho_{k(r)-1}}{2}}(Q_{k(r)-1}) $,  since $u \in
 C^{1,\a}(U^2_{2r_0})$ we have
$$\left|\frac{\partial u(\bar{x})}{\partial \nu}\right| \le
 \left|\dfrac{\partial u({x})}{\partial \nu}\right| +C|x-\bar{x}|^{\a} \le
\left|\dfrac{\partial u({x})}{\partial \nu}\right| +C\bigg( {\frac{2}{\mu} r}\bigg)^{\a}\ .$$
Integrating over $ B_{\frac{\rho_{k(r)-1}}{2}}(Q_{k(r)-1}) $, we deduce that
\begin{eqnarray}
\left|\frac{\partial u(\bar{x})}{\partial \nu}\right|^2&\le&
\frac{2}{{\omega_n}{(\frac{\rho_{k-1}}{2})}^n}\int_{B_{\frac{\rho_{k(r)-1}}{2}}\big(Q_{k(r)-1}\big)}
\bigg| \frac{\partial u(x)}{\partial \nu} \bigg|^2 \mbox{d}x + 2 C^2{\bigg(\frac{4 r^2}{{\mu}^2}\bigg)}^{\a}\le \nonumber\\
&\le&
\frac{2}{{\omega_n}{(\frac{\rho_{k-1}}{2})}^n}\int_{B_{\frac{\rho_{k(r)-1}}{2}}\big(Q_{k(r)-1}\big)}
|\nabla u(x)|^2 \mbox{d}x + 2 C^2{\bigg(\frac{4 r^2}{{\mu}^2}\bigg)}^{\a}\ .\nonumber
\end{eqnarray}
Applying the Caccioppoli inequality, we have
$$\left|\frac{\partial u(\bar{x})}{\partial \nu}\right|^2\le \frac{C}{{\big(\rho_{k-1}\big)}
^{n+2}}\int_{ B_{\rho_{k(r)-1}}(Q_{k(r)-1})} u(x)^2 \mbox{d}x +C r^{2\a} $$
and since $k$ is the smallest integer such that $d_k\le r$, then $d_{k-1}>r$, it follows

$$\left|\frac{\partial u(\bar{x})}{\partial \nu}\right|^2\le\frac{C}{\big({r \sin {\theta}^{\prime}\big)}
^{n+2}}  \int_{ B_{\rho_{k(r)-1}}(Q_{k(r)-1})} u(x)^2  \mbox{d}x + C r^{2\a}\ .$$ From $\eqref {k}$, we deduce
$$\left|\frac{\partial u(\bar{x})}{\partial \nu}\right|^2\le
\frac{C}{{r}^{n+2}} \Big\{ \big[ (\eps+E)^{1-\d}\cdot(\eps)^{\d}\big]^{{\tau}^s}
\Big\}^
{{\tau}^{k(r)-1}} +Cr^{2\a}\ .
$$
Let us define
$$ \sigma (\eps)= \big[(\eps+E)^{1-\d}\cdot(\eps)^{\d}\big]^{{\tau}^s},$$
thus the previous inequality becomes
\begin{eqnarray*}
\left|\frac{\partial u(\bar{x})}{\partial \nu}\right|^2\le
\frac{C}{{r}^{n+2}} \big\{\sigma (\eps)
\big\}^
{{\tau}^{k(r)-1}} +C r^{2\a}\ .
\end{eqnarray*}

Now, using $\eqref{r}$, we have
$$\tau^{k(r)-1}\ge \bigg(\frac{r}{d_0}\bigg)^{\nu} $$
where $\nu= -\log\big(\frac{1}{\mu}\big)\log\tau$.
We have
$$ \left|\frac{\partial u(\bar{x})}{\partial \nu}\right|\le C \bigg\{r^{-{\frac{n+2}{2}}}\Big[\sigma (\eps)\Big]^{\frac{r^{\nu}}{2}}+r^{\a} \bigg\}\ .$$
Now minimizing the function on the right hand side, with respect to $r$, with $r \in (0,\frac{r_0}{4})$, we deduce
$$\left|\frac{\partial u(\bar{x})}{\partial \nu}\right|\le C{\bigg(\log \frac{1}{\sigma(\eps)}\bigg)}
^{-\frac{2\a}{\nu+2}}\ . $$
Since this estimate holds for every $\bar{x} \in \G_{2,2r_0}$, we infer
$$ \Big\|\frac{\partial u}{\partial \nu}\Big\|_{{L^\infty}(\G_{2,2r_0})}\le C{\Big(\log{
\frac{1}{\sigma(\eps)}}\Big)}^{-\frac{2\a}{\nu+2}}$$
where $C$ is a constant depending on the \emph{a priori data} only.
Hence, solving for $\eps$, we can compute
$$\eps\ge C\exp\Big\{ {-{\Big\|\frac{\partial u}{\partial \nu}\Big\|_{L^\infty
(\G_{2,2r_0})}}}^{\!\!\!\!\!\!\!\!\!\!\!\!\!\!\!\!\!\!-\frac{\nu+2}{2\a}}\Big\}\ .
$$ Note that, recalling the a priori bound \eqref{G}, and choosing
$c=2(1-\log C G^{\gamma})$ and $\gamma=\frac{\nu+2}{2\a}$ one
trivially obtains
$$\eps\ge \exp\bigg[{-\bigg(\frac{t}{c}\bigg)^{-\gamma}}\bigg],\ \ \ \mbox{for every}\ \ t\in (0,G]\ . $$
\end{proofb}

\section{A stability result}

\begin{theorem}[Stability for a Cauchy problem]\label{stabilita'}
Let $\O$, $f_i\  i=1,2$ and $g_i$ satisfy the a priori assumptions
described above. Let $u_i \in \i$, $i=1,2$ be weak solutions of
the problem $\eqref{P}$, with $f=f_i$ and $g=g_i$
respectively and such that \eqref{E} holds for each $u_i$.\\
Moreover, let $\psi_i=u_i\big|_{\G_2},\ \ i=1,2$.
Suppose that
\begin{eqnarray*}
&& \|\psi_1 -\psi_2  \|_{L^{2}(\G_2)}\le \eps\ , \\
&& \|g_1 -g_2  \|_{L^{2}(\G_2)}\le \eps\ ,
\end{eqnarray*}
 then, for every $\rho \in (0,r_0)$
\begin{eqnarray}\label{stab}
  \| u_1 -u_2 \|_{C^1(\G_{1,\rho)}}\le \omega(\eps)
\end{eqnarray}
where $\omega$ is given by \eqref{omega} with a constant $C>0$ which depends on the \emph{a priori data} and  on $\rho$ only.
\end{theorem}

\begin{proof}
Arguing as in Theorem \ref{teo1}, we find the following estimate
\begin{eqnarray*}
&&\|u_1 -u_2 \|_{L^{\infty}(\G_{1,\frac{\rho}{2}})}\le
\tilde{\omega}(\eps)\ ,\\
&&\left\|\frac{\partial{u_1}}{\partial\nu} -\frac{\partial{u_2}}{\partial{\nu}}\right \|
_{L^{\infty}(\G_{1,\frac{\rho}{2}})}\le
\tilde{\omega}(\eps)
\end{eqnarray*}
where $\tilde{\omega}$ is a positive increasing function of the type \eqref{omega}, such that
$$\tilde
{\omega}(t)\le\tilde{C}|\log(t)|^{-\gamma_1}\ \ \mbox{for every}\ \ 0<t<1$$
where $\tilde{C}>0, 0<\gamma_1<1$ are constants depending on the \emph{a priori data} and on $\rho$ only.
By an interpolation inequality we have
$$\|\nabla_t(u_1-u_2)\|_{ L^{\infty}(\G_{1,\rho})}\le C\|u_1-u_2 \|_{ L^{\infty}(\G_{1,\frac{\rho}{2}})} ^{\beta}{\|u_1-u_2 \|_{C^{1,\a}(\G_{1,\rho})}}^{1-\beta}$$
where $\beta=\frac{\a}{\a+1}$ and $C>0$ depends on the \emph{a priori data} and on $\rho$ only,
thus by Theorem \ref{teoreg} it follows that
$$\|\nabla_t(u_1-u_2)\|_{ L^{\infty}(\G_{1,\rho})}\le C\|u_1-u_2 \|_{ L^{\infty}(\G_{1,\frac{\rho}{2}})} ^{\beta}{E}^{1-\beta}$$
where $C>0$ only depends on the \emph{a priori data} and on $\rho$ only.

It follows that for every $\eps<\eps_0$, with $\eps_0$ depending only on the \emph{a priori data}
\begin{eqnarray}
&&\|\nabla(u_1-u_2)\|_{ L^{\infty}(\G_{1,\rho})}\le\nonumber
\\
&&\le\left\| \frac{\partial{u_1}}{\partial{\nu}} -\frac{\partial{u_2}}{\partial{\nu}}\right \|
_{L^{\infty}(\G_{1,\rho})} +
\|\nabla_t(u_1-u_2)\|_{ L^{\infty}( \G_{1,\rho})}\le\nonumber
\\
&&\leq\tilde{C}\tilde{\omega}(\eps)^{\beta}
\end{eqnarray}
where $\tilde{C}>0$ depends on the \emph{a priori data} only.
Hence,
\begin{eqnarray}
\|u_1-u_2 \|_{ L^{\infty}(\G_{1,\rho})}+
\|\nabla(u_1-u_2)\|_{ L^{\infty}(\G_{1,\rho})}\le \omega(\eps)
\end{eqnarray}
where
$$\omega(t)\le C {\tilde{\omega}(t)}^{\b}\ \ \ \mbox {for every}\ \ 0< t<1\ .$$

\end{proof}
\begin{proposition}[Local monotonicity]\label{inv}
Let $u$
be a solution of \eqref{P} satisfying \eqref{E}, then there exist a point $\bar{x}
\in \G_{1,\tau}$
and a direction  $\xi \in {\mathbb{R}}^{n-1} ,\ |\xi|=1$ such that, in the representation \eqref{si} of $\G_{1}$ near $\bar{x}$, the following holds
\begin{eqnarray*}
\left|\nabla _{x'}u(x',\varphi(x'))\cdot \xi \right|\ge {\eta\left(\|g\|_{L^{\infty}(\G_{2,2r_0})} \right)} ,\
 \ x'\in U_{\bar{x}'}=\{x'=t\cdot\xi+\bar{x}',\ |t|\le\tau\}
\end{eqnarray*}
with
\begin{eqnarray}\label{tau}
\tau=\min\left\{\frac{r_0}{4} ,\frac{a \tilde{c_1}r_0}{4}, {\eta(\|g\|_{L^{\infty}(\G_{2,2r_0})})}\right\}
\end{eqnarray}
where $0<a<1, \tilde{c_1}>0$ are constants depending on the \emph{a priori data} only and $\eta$ satisfies \eqref{eta}.
\end{proposition}
\begin{proof}
Arguing as in Theorem \ref{teo1}, we can introduced a regularized distance, in the sense of Lieberman, on $S_1$ from the boundary of $\G_1$ and consequently construct connected sets $\tilde\G_{1,\rho}$ for every $\rho$, $0<\rho\le a r_0$, which satisfy
\begin{eqnarray}\label{distanza}
\G_{1,\tilde c_2h}\subset \tilde \G_{1,h}\subset \G_{1,\tilde c_1h}
\end{eqnarray}
where $0<a<1,\tilde c_2>\tilde c_1>0$ are constants depending on $M,\a$ only.

Since, by Lemma \ref{localbound}, $u \in C^{0,\a}(\G_1)$, we have that  by \eqref{distanza} it follows
$$\mathop{\rm osc}\limits_{\tilde\G_{1,\frac{\rho}{\tilde{c_1}}}} u\ge\mathop{\rm osc}\limits_{\G_{1,\frac{\tilde {c_2}\rho}{\tilde{c_1}}}} u\ge  \mathop{\rm osc}\limits_{\G_{1}} u - 2{C}E\left({\frac{\rho}{\tilde{c_1}}}\right)^{\a}{\tilde{c_2}}^{\a}\ .$$
Moreover by Theorem \ref{teo1}, we infer that
$$\mathop{\rm osc}\limits_{\tilde\G_{1,\frac{\rho}{\tilde{c_1}}}} u \ge \eta(\|g\|_{L^{\infty}(\G_{2,2r_0})})-2{C}E\left({\frac{\rho}{\tilde{c_1}}}\right)^{\a}{\tilde{c_2}}^{\a}\ .$$
Possibly replacing $c$ by a larger constant in \eqref{eta} and  taking
$$r_1=\min \left\{\eta(\|g\|_{L^{\infty}(\G_{2,2r_0})}),a\tilde{c_1}r_0, r_0\right\}$$
we have that
\begin{eqnarray}\label{osc}
\mathop{\rm osc}\limits_{\tilde\G_{1,\frac{r_1}{\tilde{c_1}}}} u \ge\eta(\|g\|_{L^{\infty}(\G_{2,2r_0})})\ .
\end{eqnarray}
Let us set, for simplicity,
$\eta=\eta\left(\|g\|_{L^{\infty}(\G_{2,2r_0})} \right)$. Since in
the a priori assumptions we have assumed that the portion $\G_1$
of the boundary is of $C^{1,\a}$ class, then we can locally
represent the restriction of $u$ (the solution to \eqref{P}) to
$\G_1$, as a function of $n-1$ variables, more precisely, for
every $x_0 \in \G_1$, up to a rigid change of coordinates, we
denote
\begin{eqnarray}\label{rapr}
w(x')= u(x',\varphi_1(x'))
 \ \ \mbox{for all } \ x\ \in \G_1 \cap B_{r_0}(x_0).
\end{eqnarray}
 By \eqref{osc}, it follows that exist two points
 $x$
and $y$ in  $ \tilde\G_{1,\frac{r_1}{\tilde{c_1}}}$,
 such that
\begin{eqnarray}\label{K}
\eta\le u(x)-u(y)\ .
\end{eqnarray}
Let us consider a continuous path $\sigma \subset \tilde\G_{1,\frac{r_1}{\tilde{c_1}}}$
joining $x$
to $y$
and let us define a sequence $\{x_i\}_{i=0,\dots,l}$ as follows
$x_0=x,\ x_i=\sigma(s_i)$ where $s_i=\max\{\mbox{$t$ },|\sigma(s)-x_i|
=\frac{r_1}{4}\}$ if $|y-x_i|>\frac{r_1}{4}$ otherwise let $i=l$ and otherwise stop the process.

The number $l$ of balls is bounded from above by $CM {\left(\frac{D}{r_1}\right)}^{n-1}$, where $C>0$ is a constant depending on $n$ only.

Let us define
$$M_i=\!\!\!\!\!\!\!\!\max_{{\overline{B_{\frac{r_1}{4}}} (x_i)
}\cap \G_1 }\!\!|\nabla_t u(x)| $$ where $\nabla_t$ denotes the
tangential gradient on $\G_1$. Let $\bar{M},\ \bar{\imath},\
\bar{x}$ be such that $\bar{x} \in B_{\frac{r_1}{4}}
(x_{\bar{i}})\cap \G_1$ and
\begin{eqnarray}\label{L}
\bar{M}=\max_{i=1,\dots,l}\left\{M_i\right\}=|\nabla_t u(\bar{x})|.
\end{eqnarray}
 By \eqref{K} and the mean value Theorem, it follows that
\begin{eqnarray*}
&\eta&\le |u(x)-u(x_1)|+\dots+|u(x_l)-u(y)|\le
\\&&\le \sum_{i=1,\dots,l} M_i\frac{r_1}{4}
\le \bar{M} C_{1}
\end{eqnarray*}
where $C_{1}>0$
is a constant depending on the \emph{a priori data} only.
Thus we have
\begin{eqnarray}\label{M}
 \bar{M}\ge \frac{\eta}{C_1}>0\ .
\end{eqnarray}
Now we use the local representation of $u$ as a function of $n-1$ variables \eqref{rapr}, within $\G_1\cap B_{\frac{r_1}{4}} (x_{\bar{i}})$.
Let us define the direction $\xi=\frac{\nabla _{x'}w}{|\nabla _{x'}w|}(\bar{x}')$.
We shall further restrict the function $w$ to the segment $t\cdot\xi+\bar{x}'$, with
$$v(t)= w(t\cdot\xi+\bar{x}')\ .$$
Now, we look for a neighborhood $U_{0}$  of $t=0$ such that

\begin{eqnarray}\label{A}
 |v'(t) |\ge \frac{\eta}{2C_1}\ \ \mbox{for every }\ \ t\in  U_{0}\ .
\end{eqnarray}

It follows that for every $|t|<\frac{r_1}{4} $
$$|v'(0)-v'(t)|\le C_2|t|^{\a}$$
where $C_2>0$ is a constant depending on the \emph{a priori data} only.

Thus  we have
$$ \bar{M}=|v'(0)|\le
 |v'(t)|  + C_2|t|^{\a}\ . $$

Hence by \eqref{M},
$$\frac{\eta}{C_1} -  C_2|t|^{\a}\le|v'(t)|\ .$$
Let us choose $t$ in such a way
$$ C_2|t|^{\a}\le\frac{\eta}{2C_1}\ . $$
Hence \eqref{A} holds with  $U_0=[-\tau,\tau]$,
where $\tau=\min\left\{\frac{r_1}{4} ,\bigg(\frac{\eta}{2C_1C_2} \bigg)^{\frac{1}{{\a}}}\right\}$.
The thesis follows, observing that $v'(t)=\frac{\partial w(x')}{\partial \xi}=\nabla_{x'}u(x',\varphi(x'))\cdot \xi$ and, possibly, by a further adjustment of the constant $c$ in \eqref{eta}.
\end{proof}

\begin{proofa}
Let $\bar{x}\in\G_{1,\tau_1},\tau_1, \xi\in \mathbb{R}^{n-1}$ be the point, the length and the direction introduced in Proposition \ref{inv}, with $u$ replaced with $u_1$. Up to a change of coordinates, we assume $\xi=e_1$.
Let $$v_i(t)=u_i(t\cdot\xi+\bar{x}, \varphi_1(t\cdot\xi+\bar{x}))\ ,\ \ \ i=1,2\ ,$$where $x=(x',\varphi_1(x'))$ is the local representation of $\G_1$ near $\bar{x}$.\\
By Proposition \ref{inv} and assumption \eqref{m}, we have that
\begin{eqnarray}\label{1}
|v_1^{\prime}(t)|\ge \eta(m)\ ,\ \ \ \  \ \mbox{for every } \ t \in U_{0}=[-\tau_1,\tau_1]\ .
\end{eqnarray}
We shall denote by $\eta_1=\eta(\|g_1\|_{L^{\infty}(\G_{2,2r_0})})$.
By the stability estimate \eqref{stab} of Theorem \ref{stabilita'}, we have that $$v'_2(t)\ge \eta_1-\omega(\eps)\ ,\ \mbox{for every}\ t\in U_0\ . $$
Thus choosing $\eps_0$ such that
$$\omega(\eps_0)\le \frac{\eta_1}{2}$$
we have
\begin{eqnarray}\label{2}
 |v_2^{\prime}(t)|\ge \frac{\eta_1}{2}\ ,\  \ \mbox{for every }  t \in U_{0}\ .
\end{eqnarray}

Thus the functions $v_i$
are invertible on $U_0$, let us denote by $V_i$ their respective images and by
\begin{eqnarray}
s^i:V_i\rightarrow U_{0}\ ,\ \ \ i=1,2\ ,
\end{eqnarray}
their inverse functions.
Let us observe that the intervals $V_1$
 and $V_2$ overlap on a sufficiently large interval $V$.
In fact, by \eqref{1} and \eqref{2} it follows that $v_i$ are monotone.Without loss of generality, let us assume they are both increasing. We have that, taken
$$a=-\frac{\tau_1}{2}\ , \ \ b=\frac{\tau_1}{2}\ ,$$
the following hold
\begin{eqnarray*}
v_i(a)<v_i(t)<v_i(b)\ ,\  \ \ \mbox{for every } t\in (a,b)\ ,\ \ i=1,2\ .
\end{eqnarray*}
Moreover, since by the Theorem \ref{stabilita'}
we have
$$ \| u_1 -u_2 \|_{L^{\infty}(\G_{1,\frac{r_1}{2}})}\le \omega(\eps) $$
then, it follows that, for $\varepsilon< \varepsilon_0$, setting $V=(v_1(a)+2\omega(\eps),v_1(b)-2\omega(\eps))$, for every $u \in  V$, there  exists $t \in(a,b) \ \mbox{such that}\ \  v_2(t)=u .$

Let us estimate from below the length of the interval $V$.
By the mean value Theorem, \eqref{1} and \eqref{tau}, it follows that
\begin{eqnarray*}
&&|v_1(a)-v_1(b)|=|v'_1(\xi)||b-a|\ge \eta_1\tau_1\ .
\end{eqnarray*}
Thus the length $\mathcal{L}$
of $V$ is bounded from below by
\begin{eqnarray*}
\mathcal{L}\ge \tau_1\eta_1-\omega(\eps)\ .
\end{eqnarray*}
 Hence, possibly adjusting the constant $c$ in the definition \eqref{eta} of $\eta$, we have that
$$\mathcal{L}\ge \eta(m)-\omega(\eps_0)\ge \frac{1}{2}\eta(m)>0\ . $$

Let us consider any value $u \in V$, then using the inverse function $s^i$,
we have
$$u=v_1(s^1(u))=v_2(s^2(u))\ .$$
Let us estimate
\begin{eqnarray*}
&&\hspace{-1cm}\!|f_1(u)-f_2(u)|=\ \nonumber\\
&&\hspace{-1cm}\!\left|\frac{\partial u_1}{\partial \nu}
(s^1(u)e_1,\varphi_1(s^1(u)e_1))\!-\!\frac{\partial u_2}{\partial \nu}
(s^2(u)e_1,\varphi_1(s^2(u)e_1))\right|\le\nonumber\\
&&\hspace{-1cm}\!\left|\frac{\partial u_1}{\partial \nu}
(s^1(u)e_1,\varphi_1(s^1(u)e_1))-\frac{\partial u_2}{\partial \nu}
(s^2(u)e_1,\varphi_1(s^2(u)e_1))\right|+ \\
&&\hspace{-1cm}\left|\frac{\partial u_2}{\partial \nu}
(s^1(u)e_1,\varphi_1(s^1(u)e_1))-\frac{\partial u_2}{\partial \nu}
(s^2(u)e_1,\varphi_1(s^2(u)e_1))\right|\hspace{+1cm}\textrm{\ \ }
\end{eqnarray*}
where $e_1=(1,0,\cdots,0)\ \in \mathbb{R}^{n-1}$.
By  Theorem
\ref{stabilita'} it follows that, for all $u \in V$,
\begin{eqnarray}\label{stima1}
&& \left|\frac{\partial u_1}{\partial \nu}
(s^1(u)e_1,\varphi_1(s^1(u)e_1))-\frac{\partial u_2}{\partial \nu}
(s^1(u)e_1,\varphi_1(s^1(u)e_1))\right|\le\omega(\eps)\ .\ \ \ \ \ \ \ \ \ \
\end{eqnarray}

 By Corollary \ref{corollario}, we infer that
\begin{eqnarray*}
&&\left|\frac{\partial u_2}{\partial \nu}
(s^1(u)e_1,\varphi_1(s^1(u)e_1))-\frac{\partial u_2}{\partial \nu}
(s^2(u)e_1,\varphi_1(s^2(u)e_1))\right|\le\\
&&\tilde{L}\left(|s^1(u)-s^2(u)|+|\varphi_1(s^1(u)e_1)-\varphi_1(s^2(u)e_1)|   \right)\le\nonumber\\
&& \tilde {L}(1+M)|s^1(u)-s^2(u)|\ .
\end{eqnarray*}
By the mean value Theorem, we find
\begin{eqnarray*}
v_2(s^2(u))=v_2(s^1(u)) + v_2'(\bar{s})(s^2(u)-s^1(u))
\end{eqnarray*}
where $\bar{s}$
is a point between $s^2(u)$ and $s^1(u)$.
Since
$$ v_2(s^2(u))=v_1(s^1(u))\ , $$
by \eqref{2} and by Theorem \ref{stabilita'}, it follows that
\begin{eqnarray}
|s^1(u)-s^2(u)|&\le& \frac{2}{\eta_1}|v_2(s^1(u))-v_1(s^1(u))|\le\nonumber\\
&\le& \frac{2}{\eta_1}\omega(\eps)\ , \ \ \ \ \ \mbox{for every}\ \ u \in V\ .\nonumber
\end{eqnarray}
Finally, we infer that
\begin{eqnarray*}
|f_1(u)-f_2(u)|\le\omega(\eps)\ ,\ \ \ \ \mbox{for every} \ \ u \in\  V\ ,
\end{eqnarray*}
possibly by a further adjustment of the constant $C$ in \eqref{omega}.
\end{proofa}

\end{document}